\documentclass{article}
\usepackage{hchang}
\usepackage{hyperref}
\usepackage{comment}
\usepackage{tcolorbox}
\usepackage{color}
\usepackage{graphicx} %
  \setkeys{Gin}{width=\linewidth,totalheight=\textheight,keepaspectratio}
  \graphicspath{{graphics/}} %
\usepackage{booktabs} %
\usepackage{units}    %
\usepackage{multicol} %
\usepackage{fancyvrb} %
  \fvset{fontsize=\normalsize}%
\usepackage{cleveref}
\usepackage{rotating}
\usepackage{transparent}
\usepackage{tikz}
\usetikzlibrary{arrows.meta}
\usepackage{fullpage}
\usepackage{soul}
\usepackage{thmtools}
\usepackage{thm-restate}
\usepackage{algorithm}
\usepackage{algorithmicx}
\usepackage{algpseudocode}

\usepackage{subcaption}

\newcommand{\blue}[1]{{\color{blue}#1}}

\newtheorem{lemma}{Lemma}
\newtheorem{theorem}{Theorem}
\newtheorem{claim}{Claim}
\newtheorem{definition}{Definition}

\newtheorem{corollary}{Corollary}

\newtheorem{conjecture}{Conjecture}

\newcommand{\pr}{\mathrm{Pr}}
\newcommand{\draw}{\in_{R}}

\newcommand{\cX}{\mathcal{X}}

\newcommand{\real}{\mathbb{R}}

\newcommand{\defn}[1]{\EMPH{#1}}

\newcommand{\bbE}{\mathbb{E}}
\newcommand{\cK}{\mathcal{K}}

\crefname{problem}{problem}{problems}
\Crefname{problem}{Problem}{Problems}
\crefformat{problem}{problem~(#2#1#3)}
\Crefformat{problem}{Problem~(#2#1#3)}
\usepackage{xcolor}
\newcommand{\lightgray}[1]{\textcolor{gray}{#1}}

\renewcommand{\leq}{\leqslant}
\renewcommand{\le}{\leqslant}
\renewcommand{\geq}{\geqslant}
\renewcommand{\ge}{\geqslant}
\newcommand{\bdry}{\partial}

\algtext*{EndIf}
\algtext*{EndFor}
\algtext*{EndWhile}
\algtext*{EndProcedure}

\algrenewcommand\algorithmicrequire{\textbf{Input:}}
\algrenewcommand\algorithmicensure{\textbf{Output:}}
\makeatletter
\newcommand{\algmultiline}[1]{\parbox[t]{\dimexpr\linewidth-\ALG@thistlm}{#1\strut}}
\makeatother

\newcommand{\dist}{\ensuremath{d}}
\nolinenumbers

\title{$K_h$-Minor-Free Graphs Have Balanced Separators of Size $O(h\sqrt{n})$}

\author{
Jonathan Conroy\thanks{Dartmouth College. \texttt{jonathan.conroy.gr@dartmouth.edu}. Supported by the NSF CAREER award CCF-2443017.}
\and
Hung Le\thanks{University of Massachusetts Amherst.  \texttt{hungle@cs.umass.edu}. Supported by the NSF CAREER award CCF-2237288 and the NSF grant CCF-2517033.}
\and
Jack Spalding-Jamieson\thanks{Independent. \texttt{jacksj@uwaterloo.ca}}
}
\date{}

\begin{document}

\maketitle

\begin{abstract}

We show that every $K_h$-minor-free graph with $n$ vertices has a balanced separator of size $O(h\sqrt{n})$ for any parameter $h\in\mathbb{N}$, resolving the conjecture posed by Alon, Seymour, and Thomas~\cite{AST90} in 1990. The proof is elementary, building on the combinatorial framework by Le~\cite{Le26}. A key new idea of the proof is the existence of a \EMPH{short hitting walk} on vertex expanders: given $k$ large  (not necessarily disjoint) subsets of vertices, there is a walk of length $O(k/\psi(G))$ hitting every subset, where $\psi(G)$ is the vertex expansion. This technical result substantially improves on an analogous lemma by Alon, Seymour, and Thomas~\cite{AST90} in an important parameter setting, where the previous lower bound~\cite{Kotlov00} does not apply.

\end{abstract}

\paragraph{AI Disclosure.} The key technical ideas of this paper
are due to GPT-6 Astra, with some simplification due to Claude Opus 5.5.
A Lean formalization\footnote{\url{https://github.com/jacketsj/formalized-minor-free-separators}}
of our main result was produced using a combination of GPT-6 Astra and Claude Opus 5.5, and its main definitions and statements have been checked over by the human authors.
GPT-6.1 Sol was used to produce some figures from hand-drawn sketches.
The human authors claim no role in developing the technical ideas. Our main contribution is verifying correctness, contextualizing the results, completely rewriting the paper to present the key new ideas in the proofs, and attributing existing ideas to prior work (to the best of our knowledge). In particular, we have simplified the proof for understandability,
and obtained a slightly worse constant than the formalization in doing so.

\section{Introduction}

The separator theorem for planar graphs by Lipton and Tarjan~\cite{LT79} in 1979 is a cornerstone of algorithmic graph theory. Given a graph $G = (V,E)$, a vertex subset $S\subset V$ is a \EMPH{balanced separator} if every connected component of $G\setminus S$ has size at most $3|V|/4$. The Lipton-Tarjan theorem says that any $n$-vertex planar graph has a balanced (vertex) separator of size $O(\sqrt{n})$. This bound is asymptotically optimal: any balanced separator of the $\sqrt{n}\times \sqrt{n}$-planar grid must have size $\Omega(\sqrt{n})$. A similar result was soon shown for bounded genus graphs~\cite{GHT84}, which is also asymptotically tight. A decade later, Alon, Seymour, and Thomas~\cite{AST90} (AST) generalized the planar separator theorem to $K_h$-minor-free graphs. They showed that any $n$-vertex $K_h$-minor-free graph has a balanced separator of size $O(h^{3/2}\sqrt{n})$. In particular, for a constant $h$, the AST separator has size $O(\sqrt{n})$; this applies to planar graphs, which are \(K_5\)-minor-free. In the regime of constant $h\geq 5$, their separator bound is also asymptotically optimal, witnessed by the  $\sqrt{n}\times \sqrt{n}$-planar grid.

For non-constant $h$, Alon, Seymour, and Thomas~\cite{AST90} conjectured that the bound is not optimal. They pointed out that for $3$-regular expander graphs with $n$ vertices, $h = \Theta(n^{1/2})$, but the bound $O(h^{3/2}\sqrt{n}) = O(n^{1/4}n)$, while a (trivial) balanced separator has size $n$. This excess of $\sqrt{h}$ is their basis for the following conjecture.

\begin{conjecture}[Alon-Seymour-Thomas Conjecture~\cite{AST90}]\label{conj:separator} Any $K_h$-minor-free graph with $n$ vertices has a balanced vertex separator of size $O(h\sqrt{n})$ for any $h\geq 3$.
\end{conjecture}

\Cref{conj:separator} has received significant attention over the last few decades; see \Cref{tab:history} for an overview. In 1994, Plotkin,  Rao, and Smith~\cite{PRS94} constructed a balanced separator of size $O(h\sqrt{n \log n})$. For large $h$, this improves on the AST bound, but for constant $h$, their separator size is not optimal. Nearly two decades later, in 2010, Kawarabayashi and Reed~\cite{KR10}  announced a separator result based on the Robertson–Seymour structure theorem; their precise bound is later clarified\footnote{See \cite{Le26,bonnet2026separator} for a more detailed discussion.}  to be $O(h\sqrt n+f(h))$, where $f(h)$ is a very fast-growing function. It is possible that the recently improved structure theorem for minor-free graphs~\cite{GSW25} may reduce $f(h)$ to a (very large) polynomial in $h$, but even then it does not settle \Cref{conj:separator}, and the proof would not be elementary.
Around the same time,  Biswal, Lee, and Rao introduced a new technique to give another balanced separator of size $O(h^3\sqrt{\log h}\cdot\sqrt{n})$. Although their bound is worse than the AST bound, their technique brought multicommodity flow-cut duality into the picture. Recent works~\cite{SpJ25,Le26} exploit the optimal padded decomposition~\cite{CF25} to push the flow-cut duality technique much further and obtain a balanced separator of size $O(h\log h\cdot\sqrt{n})$~\cite{Le26}.  This bound matches the barrier of the standard flow-cut approach, since the flow-cut gap can be $\Omega(\log h)$
(e.g. for $3$-regular expanders).

\begin{table}[ht]
\centering
\begin{tabular}{llcp{6.2cm}}
\toprule
Separator size & Reference & Year & Notes \\
\midrule
$O(h^{3/2}\sqrt{n})$ & \cite{AST90} & 1990 & \\
$O(h\sqrt{n\log n})$ & \cite{PRS94} & 1994 & \\
$O(h^3\sqrt{\log h}\cdot\sqrt{n})$ & \cite{BLR10} & 2008 & flow-cut duality \\
$O(h\sqrt{n} + f(h))$ & \cite{KR10} & 2010 & $f(h$) is a fast growing function \\
$O(h\log h\log\log h\cdot\sqrt{n})$ & \cite{SpJ25} & 2025 & \\
$O(h\log h\cdot\sqrt{n})$ & \cite[Thm.~1]{Le26} & 2026 & optimal flow-cut duality \\
$O(h\sqrt{\log h}\cdot\sqrt{n})$ & \cite[Thm.~2]{Le26} & 2026 & bypasses flow-cut duality \\
$O(h\sqrt{n})$ & \multicolumn{2}{c}{This paper} & Asymptotically optimal~\cite{AST90} \\
\bottomrule
\end{tabular}
\caption{Bounds for balanced separators of $n$-vertex $K_h$-minor-free graphs.}
\label{tab:history}
\end{table}

Recently, Le~\cite{Le26} bypassed the flow-cut duality barrier and obtained a balanced separator of size  $O(h\sqrt{n \log h})$. The main idea there is to incorporate the optimal padded decomposition by Conroy and Filtser~\cite{CF25} directly into the combinatorial frameworks by Alon, Seymour, and Thomas~\cite{AST90}  and Plotkin,  Rao, and Smith~\cite{PRS94}. The $\sqrt{\log h}$ factor is essentially $\sqrt{\beta_h}$ where $\beta_h = \Theta(\log h)$ is the ``padding parameter'' of $K_h$-minor-free graphs~\cite{CF25}. The bound $\log h$ is also optimal for $\beta_h$. As noted in~\cite{Le26}, ``Erasing $\sqrt{\log h}$ entirely seems to require substantially new ideas.''

Our main result is a complete resolution of \Cref{conj:separator}. 

\begin{restatable}{theorem}{SepMain}\label{thm:sep-main}
Any $n$-vertex $K_h$-minor-free graph $G$ has a $3/4$-balanced separator of size at most $16h\sqrt{n}$.
\end{restatable}

In \Cref{subec:overview}, we give a detailed overview of the proof. Here we highlight two new ideas underlying the proof. 

\paragraph{The first idea} is a bound on the length of a walk intersecting $k$ large subsets of vertices in a vertex expander. To formally state the bound, we recall the definition of a vertex expander. Given a connected graph $G = (V,E)$ and a subset $S\subseteq V$ of vertices, we define the \defn{(outer) boundary} $\text{\EMPH{$\bdry S$}} \coloneqq \{v\in V\setminus S:\exists u\in S,uv\in E\}$, which is the set of vertices in $V\setminus S$ that have at least one neighbor in $S$. We say that $G$ is a \EMPH{$\psi$-vertex expander} for a parameter $\psi > 0$ if, for every $S\subseteq V$ such that $|S| \leq |V|/2$, $|\bdry S|\geq \psi |S|$. That is, every vertex subset of size at most $|V|/2$ has a sufficiently large external neighbor set in $G$. We denote by $\text{\EMPH{$\psi(G)$}} = \max\{\psi: G \text{ is a $\psi$-vertex expander}\}$.

The central new ingredient of the proof of \Cref{thm:sep-main} is the following theorem, which gives a short walk that intersects every set in a given collection of large subsets. Note that these subsets may not be vertex-disjoint.

\begin{restatable}[Short Hitting Walk on Expanders]{theorem}{HittingThm}
\label{thm:hitting} Let $G=(V,E)$ be a connected graph on $n$ vertices. Let $X_1,\dots,X_k\subset V$ be $k$ subsets of vertices, each of size at least $n/k$. Then $G$ has a walk $W$ of length
less than $13k/\psi(G)$ that intersects every subset: $W\cap X_i\not=\varnothing$ for every $i \in [k]$.
\end{restatable}

To place \Cref{thm:hitting} in context, we compare it with the following lemma of Alon, Seymour, and Thomas~\cite{AST90}.

\begin{lemma}[Lemma 2.1~\cite{AST90}]\label{lm:walk-AST} Let $G = (V,E)$ be an $n$-vertex graph and $\ell > 0$ be a parameter. Let $X_1,\ldots, X_k$ be $k$ subsets of vertices of $G$. Then there exists either:
\begin{enumerate}
    \item a tree $T$ in $G$ of at most $\ell k$ vertices such that $T\cap X_i \not=\varnothing$ for every $i\in [k]$, or
    \item  a subset $Z\subseteq V$ of size $|Z| < n/\ell$ such that for every connected component $C$ of $G\setminus Z$, $C\cap X_i = \emptyset$ for some $i\in [k]$.
\end{enumerate}
\end{lemma}

 Alon, Seymour, and Thomas~\cite{AST90} observed that improving the trade-off between the tree size $h\ell$ and the separator size $n/\ell$ would improve their balanced-separator bound (which is $O(h^{3/2}\sqrt{n})$). However, Kotlov~\cite{Kotlov00} later showed that the tradeoff is tight when $\ell$ is constant, and $G$ is the hypercube.  Subsequent papers either avoid the AST framework entirely~\cite{KR10,BLR10,SpJ25} or exploit $K_h$-minor-freeness to get a better tradeoff~\cite{Le26}. (Herein, for simplicity, we \underline{assume that $k = h$} since this is the most interesting regime.) 
 Specifically, Le~\cite{Le26} uses the optimal padded decomposition for $K_h$-minor-free graphs ~\cite{CF25} to construct a tree $T$ of size $O(\ell \log h)$ (in a very specific setting), improving over the bound $O(h\ell)$ in \Cref{lm:walk-AST}. The $\log h$ factor is essentially the padding parameter, which is best possible as noted above.

The walk-length bound in \Cref{thm:hitting} improves the bound on $|V(T)|$ in \Cref{lm:walk-AST} to $O(\ell)$ when $\ell = c\cdot h/\psi(G)$ for a sufficiently large constant $c$. %
To see this, assume $|X_i| = n/h$ for every $i\in [k]$, and let $Z\subseteq V$ be any set with
$|Z| < n/\ell = \psi(G)\cdot n/(c\cdot h)$.
Every union $S$ of components of $G\setminus Z$ with $|S|\le n/2$ has $\bdry S\subseteq Z$,
so $|S|\le |Z|/\psi(G) < n/(c\cdot h)$ by the definition of a vertex expander.
Take such an $S$ of maximum size.
By maximality, every other component $C$ of $G\setminus Z$ has $|C| > n/2-|S|>n/2-n/(c\cdot h)$,
which is greater than $n/(c\cdot h)>|S|$ for $c\geq 2$,
and therefore $|C| > n/2$.
So $G\setminus Z$ consists of $S$ and a single component $C$,
which intersects every $X_i$ because
$|V\setminus C| = |S|+|Z| < n/(c\cdot h) + \psi(G)\cdot n/(c\cdot h) = (1+\psi(G))\cdot n/(c\cdot h) \leq n/h = |X_i|$ for any $c\geq 3\geq1+\psi(G)$.
Hence $Z$ does not satisfy case 2 of \Cref{lm:walk-AST}, which requires every component of $G\setminus Z$ to miss some $X_i$.
Hence, case 2 cannot occur,
and case 1 gives a tree $T$ with at most
\begin{equation*}
    h\ell  = O(h^2/\psi(G))
\end{equation*}
vertices that intersects every $X_i$.
However, \Cref{thm:hitting} gives a tree --- more precisely, a spanning tree of the subgraph $G[V(W)]$ where $W$ is the walk --- of only $O(h/\psi(G))$ vertices that intersects every $X_i$, a factor of $h$ improvement!
Note that \Cref{thm:hitting} does not exploit $K_h$-minor-freeness, and yet still gives a factor of $h$ improvement over \Cref{lm:walk-AST}.

At this point, readers might ask: does \Cref{thm:hitting} contradict Kotlov's~\cite {Kotlov00} lower bound, which shows the tightness of \Cref{lm:walk-AST}? One observation is that \Cref{thm:hitting} requires $|X_i|\geq n/k$, whereas \Cref{lm:walk-AST} imposes no such requirement. However, the lower bound instance of  Kotlov~\cite {Kotlov00} has $|X_i| = n/2$, so assuming large subsets is not the source of the improvement. The key fact here is that the parameter setting $\ell = c\cdot k/\psi(G)$ in the application of \Cref{thm:hitting} is \EMPH{outside the lower bound regime of Kotlov}. Specifically,  Kotlov's lower bound applies to the case $\ell = 2$, $k \approx 2\log(n)$ and $G$ is the (discrete) hypercube (in $\{0,1\}^{\log n}$). For the hypercube\footnote{For the hypercube, \Cref{thm:hitting} is asymptotically tight when $k=4$.}, $\psi(G) = \Theta(1/\sqrt{\log n})$~\cite{Harper66}, and hence our setting of $\ell$ in \Cref{thm:hitting}  is $\ell = c\cdot k/\psi(G) = \Theta(\log^{3/2}(n))$, which is much larger than $2$. In summary, there is no contradiction.

\paragraph{The second idea} is to adapt the algorithm by Le~\cite{Le26}, which itself is a combination of two different algorithms~\cite{AST90,PRS94}, so that \Cref{thm:hitting} can be applied. For readers who are familiar with  Le~\cite{Le26} (as well as ~\cite{AST90,PRS94}), in a nutshell, the algorithm is iterative, and every iteration essentially has two major steps: it finds either (i)  a small set of vertices $Z$ of the current subgraph $H\subseteq G$ that cuts off a large subset $A$ (relative to $|Z|$), or (ii) a tree  $T\subseteq H$ of small size that intersects \EMPH{every neighborhood} of the current branch sets in the clique minor model. Here in this paper, both steps are reframed in terms of vertex expanders. Specifically, case (i)  can be implemented as cutting off the smaller side of a sparse cut if $H$ is not a vertex expander, while in the (complementary) case (ii), $H$ is a vertex expander. In this case, we split into two smaller cases: one of which extends the neighborhood of each branch set, so its size is at least $n/h$, then applies \Cref{thm:hitting} to find the tree $T$. The other smaller case is handled similarly to prior work~\cite {AST90,Le26}.

\paragraph{On notation.} Graphs in this paper are unweighted and connected, unless explicitly noted otherwise.  For any walk $W$ in $G$, we denote by $|W|$ the length of the walk, which is the number of edges in $W$ (counted with multiplicity).  Let $\pi_G(u,v)$ denote a shortest path between $u$ and $v$ in a graph $G$. We use $d_G(u,v)$ to denote the distance between $u$ and $v$, which is also $|\pi_G(u,v)|$. For any two walks $P$ and $Q$ of $G$ that share the same endpoint, we denote by $P\circ Q$ the walk obtained by concatenating $P$ and $Q$ at the shared endpoint. We use $[n]$ to denote the set $\{1,2,\ldots, n\}$.

\subsection{Related work}

An orthogonal line of work on balanced separators of $K_h$-minor-free graphs focuses on fast, ideally linear-time algorithms~\cite{AST90,PRS94,RW09,KR10,WulffNilsen11,wulffnilsen14,bonnet2026separator}. Here, we assume that the minor size $h$ is a fixed constant. In this regime, Bonnet, Korhonen, Le, Li, and Masařík~\cite{bonnet2026separator} gave a linear-time algorithm, which concludes decades of effort on the running-time question. The dependency on $h$ in these works could be exponential~\cite{RW09} or polynomial with a high degree~\cite{bonnet2026separator}.

Another relevant line of work focuses on constructing balanced separators for more general graph classes, such as graphs with polynomial expansion~\cite{PRS94,DN16,ER18}, $H$-induced-minor-free graphs~\cite{KL24}, or fat-minor-free graphs~\cite{BLPP26}. In these papers, the separator size, especially its dependency on the minor size (in their respective definitions of graph minors), is not optimal. It would be interesting to see if the short hitting walk in \Cref{thm:hitting} can be applied to improve the separators in these more general settings.

\subsection{Proof overview} \label{subec:overview}

A \EMPH{$K_t$ minor model} of a graph $G$ is a collection $\cK = \{C_1,C_2,\ldots, C_t\}$ of vertex-disjoint connected subgraphs of $G$, such that for any two different sets $C_i, C_j$ with $1\leq i \not=j \leq t$ there exists an edge $(u,v)\in E$ where $u\in C_i$ and $v\in C_j$.  Each subgraph in $\cK$ is called a \EMPH{branch set}. We use $V(\cK)$ to denote the set of vertices $\bigcup_{C \in \cK} C$.

\paragraph{The algorithm.} We describe an algorithm to find an $O(h \sqrt n)$-size separator. (For simplicity, we are concerned only with the existence of the separator and not an efficient construction. The algorithm as we present it does not run in polynomial time.) In describing the algorithm, we assume that the walk-length bound in \Cref{thm:hitting} holds. The algorithm is basically a reduction to the case where we can apply the walk-length bound. At the end, we give an overview of the proof of \Cref{thm:hitting}.

Just as in \cite{AST90, PRS94, Le26}, our algorithm proceeds by iteratively sweeping a short separator through $G$ until it finds a balanced separator. It maintains an induced subgraph \EMPH{$H$} of $G$ which contains at least $3n/4$ vertices; the \EMPH{current separator} is the boundary of $H$, denoted \EMPH{$\bdry H$}. One should think of $H$ as the largest connected component of $G$ after removing the ``current separator.'' We maintain the invariant that the current separator, i.e.\ the boundary of $H$, always has size $O(h \sqrt n)$. To this end, the algorithm maintains two additional objects, whose union contains the boundary of $H$: (1) a minor model \EMPH{$\cK$} of a clique in $G$, such that each branch set $C \in \cK$ has $|\bdry H \cap C| \le O(\sqrt n)$
and (2) an additional subset \EMPH{$X$} of vertices, which we separately guarantee has size $O(h \sqrt n)$. The invariant that $\bdry H$ lies in the union $X \cup V(\cK)$ implies that $|\bdry H| \le O(h \sqrt n)$.

In each iteration of the algorithm, we add some vertices to $X \cup V(\cK)$, then set $H$ to be the largest connected component of $H \setminus (X \cup V(\cK))$. When the size of $H$ drops below $3n/4$, we return the boundary of $H$ as our balanced separator. In each iteration, there are three cases that determine which vertices to add to $X \cup V(\cK)$. The first two are simple cases. In the only remaining case, we can exploit the expander walk theorem.

\begin{itemize}
    \item \textbf{Case 1: $H$ has small expansion.}

    Case 1 occurs if $H$ is \emph{not} an $h/\sqrt n$-expander. In this case, there is a set of vertices $A$ whose boundary $Z$ contains at most $h/\sqrt n |A|$. We add $Z$ to $X$, and continue to the next iteration. A charging argument\footnote{Each vertex in $A$ receives $h/\sqrt n$ charge from the vertices of $Z$, and is never charged again. Over $n$ vertices, there is at most $h \sqrt n$ charge.} shows that $X$ never grows beyond $O(h \sqrt n)$.

\item \textbf{Case 2: Some branch set $C \in \cK$ can be enlarged while preserving $|\bdry H \cap C| \le O(\sqrt n)$. }

For Case 2, we check if there is some branch set $C$ and radius $r$ such that we can enlarge $C$ by growing a radius-$r$ ball around $C$ into $H$. Such an enlargement is valid so long as $C$ still satisfies $ |\bdry H \cap C| \le O(\sqrt n)$, which is equivalent to asking if there are at most $O(\sqrt n)$ vertices at distance $r$ from $C$ (with respect to $H$). If there is such a branch set $C$ and small-boundary ball, then we are in Case 2: we grow $C$ and update $\cK$.

\item \textbf{Case 3. Otherwise, we can apply the walk theorem.}

If Case 1 does not apply, we know $H$ is an $h/\sqrt n$-expander. If Case 2 does not apply, we can show that for any branch set $C \in \cK$, the $\sqrt n/ h$-radius ball
around $C$ contains at least $n/h$ vertices in $H$. In this case, the hitting walk theorem applies, and we can find a walk $W$ in $H$ with $O(\sqrt n)$ vertices that travels within distance $\sqrt n/h$ of every $C \in \cK$. We extend the walk $W$ so that it is adjacent to every $C \in \cK$. We have that $|W| = O(\sqrt n)$. The vertices in $W$ become a new branch set that we add to $\cK$.
\end{itemize}
In each of these three cases, we shrink $H$ while maintaining the invariant that $|\bdry H| \le O(h \sqrt n)$. When $H$ becomes small enough, we terminate and return $\bdry H$.

\paragraph{The existence of a short hitting walk in \Cref{thm:hitting}.} It is well known that a $\psi(G)$-vertex-expander has diameter $O(\log n/\psi(G))$, which is also tight.  Thus, we can construct a walk in \Cref{thm:hitting} of length $O(k\log n/\psi(G))$ simply by picking one vertex per subset and connecting them into a walk using shortest paths. This bound is worse than \Cref{thm:hitting} by a factor of $\log n$ and does not exploit the fact that the subsets are large (of size at least $n/k$). 

Getting rid of this $\log n$-factor is more complicated.
The very high-level idea of the proof is to construct a randomized walk \EMPH{$W$} that performs exactly $k$ large steps, using shortest-paths between chosen vertices,
while selecting exactly one of the sets $X_1,\dots,X_k$ in each step.
We design the random walk to minimize the \emph{expected length} of visiting a new subset.
For simplicity of the exposition in this overview, we assume that the subsets $\{X_1,\ldots, X_k\}$ are \EMPH{pairwise vertex disjoint}, and each has size exactly \EMPH{$r \coloneqq \ceil{n/k}$}.

In this special case, the random walk maintains the following invariant: After the $t$-th step of the walk, we have selected a uniformly random $t$-element subset $I \subseteq [k]$, and the walk is guaranteed to hit all sets $X_i$ with $i \in I$. Moreover, conditioned on $I = J$ (for any fixed $t$-element subset $J \subseteq [k]$), the last vertex on the walk is distributed uniformly at random among the vertices in the selected sets $\set{v \in X_i : i \in J}$.

How should we design the random walk to maintain this invariant? Fix some $J \subseteq [k]$ and condition on the fact that our selected subset is $I = J$. We know that the last vertex $v_t$ is uniformly distributed among $\set{v \in X_i : i \in J}$. Call this distribution \EMPH{$p_J$}. To maintain the invariant, we would like to select a uniformly random index $i \in [k] \setminus J$ and a uniformly random vertex $v_{t+1} \in X_i$, and then extend the walk $W$ along a shortest path from $v_t$ to $v_{t+1}$. Let \EMPH{$q_J$} denote the distribution from which $v_{t+1}$ is sampled. It is not difficult to see that if we sample $v_{t+1}$ according to $q_J$, we maintain the desired invariant on the random walk.

The simplest way to sample $v_{t+1}$ from $q_J$ is to perform this sampling independently of the current vertex $v_t$ (which itself was sampled from $p_J$). However, this turns out to be a bad idea: the expected distance $\E[d(v_t, v_{t+1})]$ is too large, and we recover the same $O(k\log n/\psi(G))$ bound. Instead, we can sample from $q_J$ in a way that is \EMPH{correlated} with $p_J$. We exploit the fact that $G$ is an expander to show\footnote{Specifically, minimizing the expected distance while sampling jointly from $p_J$ and $q_J$ is equivalent to bounding the \emph{Wasserstein distance} between $p_J$ and $q_J$. We bound the Wasserstein distance on expanders by using a standard dual form of Wasserstein distance from optimal transport literature,
and a known relationship of that dual form to vertex expansion.
} that $v_{t+1}$ can be sampled from $q_J$ such that 
\[\E[d(v_t, v_{t+1})] \le \frac{4k}{\psi(G) \sqrt{t(k-t)}}.\]
By linearity of expectation and a routine calculation, the expected length of the walk can then be bounded by
\begin{equation}
      \bbE[|W|] \leq  \frac{4k}{\psi(G)}
\sum_{j=1}^{k-1}\frac1{\sqrt{j(k-j)}}  \leq \frac{4k}{\psi(G)}
\int_0^k\frac{dt}{\sqrt{t(k-t)}} =  \frac{4\pi k}{\psi(G)} < \frac{13k}{\psi(G)}
\end{equation}
as claimed in \Cref{thm:hitting}.

In this overview, we used the pairwise disjointness assumption on the subsets to simplify the definitions of the two distributions $p_J$ and $q_J$. In the actual proof, we slightly adjust the definitions of $p_J$ and $q_J$ to account for overlaps between subsets. Consequently, the proof that $\E[d(v_t, v_{t+1})] \le \frac{4k}{\psi(G) \sqrt{t(k-t)}}$ is somewhat more complicated than it would be in the disjoint-subset case.

\section{Short Hitting Walks on Expanders}

In this section, we prove \Cref{thm:hitting}, restated below.

\HittingThm*

We prove \Cref{thm:hitting} by constructing a randomized walk that can be viewed as a randomized greedy process in $k-1$ steps. Let $\cX = \{X_1,X_2,\ldots, X_k\}$. By discarding vertices from the subsets, we can assume that every subset $X\in \cX$ has size exactly \EMPH{$r \coloneqq \ceil{n/k}$}. For a subset $I \subseteq [k]$, we denote by $\cX_I = \{X_i: i\in I\}$ the set of subsets whose indices are in $I$.

\paragraph{The algorithm.} We start the walk $W$ at a random vertex $v_1\in X_i$, where $X_i$ is sampled uniformly at random from $\{X_1,\ldots, X_k\}$. At the beginning of step $t\in [k-1]$, let $I_t$ be the set of indices selected so far. We will maintain that $|I_t| = t$, and that the walk $W$ intersects every subset in $\cX_{I_t}$. Let $\Bar{I}_t = [k]\setminus I_t$, which is the set of indices not yet selected.

 \usetikzlibrary{backgrounds}

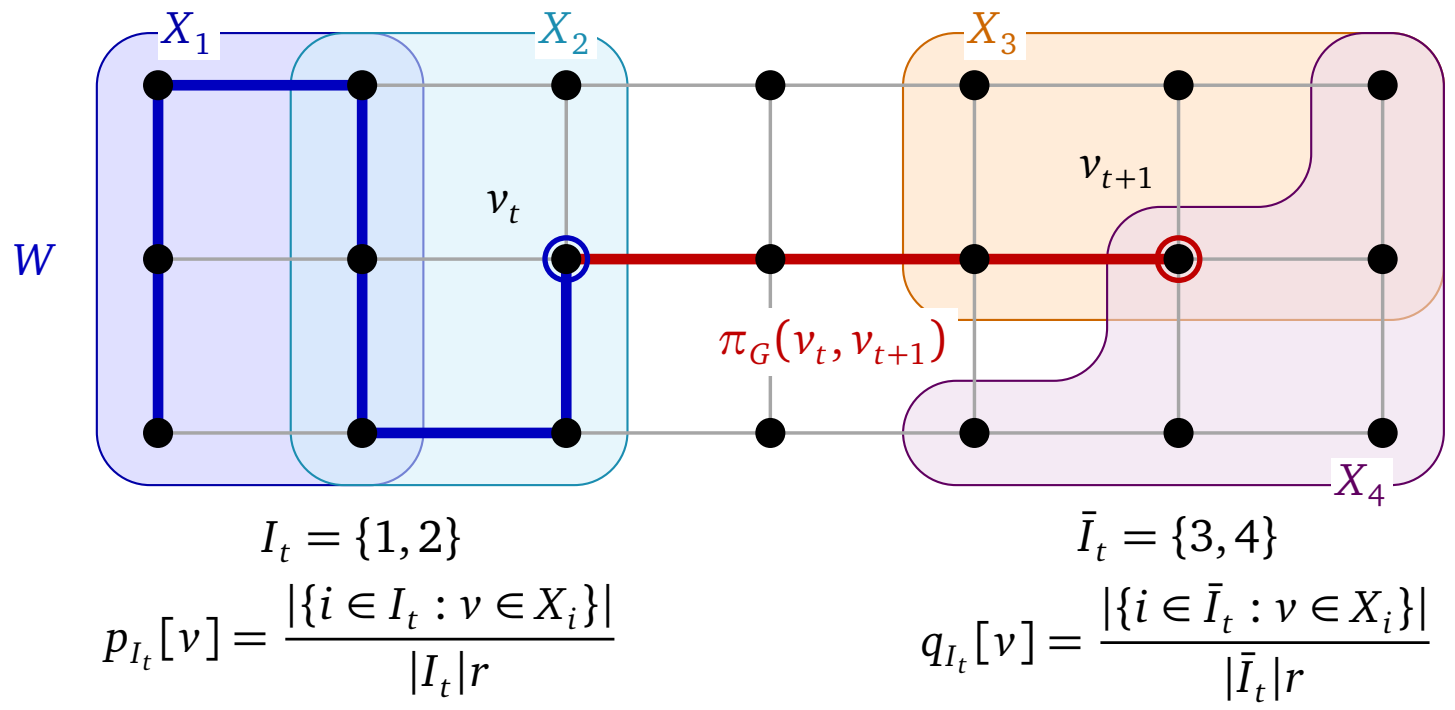
\begin{figure}[t]
\centering
\begin{tikzpicture}[
    x=1.35cm,y=1.15cm,
    vertex/.style={
        circle,fill=black,inner sep=2pt
    },
    edge/.style={
        draw=black!35,line width=.65pt
    },
    oldwalk/.style={
        draw=blue!75!black,line width=2pt
    },
    extension/.style={
        draw=red!75!black,line width=2pt
    },
    every node/.style={font=\small}
]

\foreach \row in {0,...,2}{
    \foreach \col in {0,...,6}{
        \coordinate (v\col\row) at (\col,\row);
    }
}

\begin{scope}[on background layer]
    \filldraw[
        rounded corners=10pt,
        draw=blue!65!black,fill=blue!12
    ]
        (-.3,-.3) rectangle (1.3,2.3);

    \filldraw[
        rounded corners=10pt,
        draw=cyan!65!black,fill=cyan!15,
        fill opacity=.55
    ]
        (.65,-.3) rectangle (2.3,2.3);

    \filldraw[
        rounded corners=10pt,
        draw=orange!80!black,fill=orange!15
    ]
        (3.65,.65) rectangle (6.3,2.3);

    \filldraw[
        rounded corners=10pt,
        draw=violet!75!black,fill=violet!15,
        fill opacity=.55
    ]
        (3.65,-.3) -- (6.3,-.3) -- (6.3,2.3) --
        (5.65,2.3) -- (5.65,1.3) -- (4.65,1.3) --
        (4.65,.3) -- (3.65,.3) -- cycle;
\end{scope}

\foreach \row in {0,...,2}{
    \foreach \col [evaluate=\col as \next using int(\col+1)]
        in {0,...,5}{
        \draw[edge] (v\col\row) -- (v\next\row);
    }
}
\foreach \col in {0,...,6}{
    \foreach \row [evaluate=\row as \next using int(\row+1)]
        in {0,1}{
        \draw[edge] (v\col\row) -- (v\col\next);
    }
}

\draw[oldwalk]
    (v00) -- (v01) -- (v02) -- (v12)
    -- (v11) -- (v10) -- (v20) -- (v21);

\draw[extension]
    (v21) -- (v31) -- (v41) -- (v51);

\foreach \row in {0,...,2}{
    \foreach \col in {0,...,6}{
        \node[vertex] at (v\col\row) {};
    }
}

\node[
    circle,draw=blue!75!black,
    line width=1pt,minimum size=8pt,inner sep=0pt
] at (v21) {};
\node[
    circle,draw=red!75!black,
    line width=1pt,minimum size=8pt,inner sep=0pt
] at (v51) {};

\node at (1.7,1.3) {$v_t$};
\node at (4.7,1.5) {$v_{t+1}$};

\node[text=blue!65!black,fill=white,inner sep=1pt]
    at (.15,2.3) {$X_1$};
\node[text=cyan!65!black,fill=white,inner sep=1pt]
    at (2,2.3) {$X_2$};
\node[text=orange!80!black,fill=white,inner sep=1pt]
    at (4.1,2.3) {$X_3$};
\node[text=violet!75!black,fill=white,inner sep=1pt]
    at (5.9,-.3) {$X_4$};

\node[text=blue!75!black] at (-0.6,1) {$W$};
\node[
    text=red!75!black,below=9pt,
    fill=white,inner sep=2pt
] at (3.3,1) {$\pi_G(v_t,v_{t+1})$};

\node[align=center] at (1,-1.0) {
    $I_t=\{1,2\}$\\[3pt]
    $\displaystyle
    p_{I_t}[v]=\frac{|\{i\in I_t:v\in X_i\}|}{|I_t|r}$
};
\node[align=center] at (5,-1) {
    $\bar I_t=\{3,4\}$\\[3pt]
    $\displaystyle
    q_{I_t}[v]=\frac{|\{i\in\bar I_t:v\in X_i\}|}
                    {|\bar I_t|r}$
};

\end{tikzpicture}
\caption{
One iteration of the walk construction on a graph.
Each shaded subset contains $r=6$ vertices. The red path extends $W$ from $v_t$ to $v_{t+1}$. Vertices belonging to two subsets in the same family receive
twice the probability mass of vertices belonging to only one.
}
\label{fig:wasserstein-walk}
\end{figure}

The walk $W$ currently ends at a vertex we denote as $v_t$; we can show (in \Cref{lemma:distribution}) that the distribution of $v_t$ conditional on $I_t = J$ (for some subset $J \subset [k]$) is:
\begin{equation}\label{eq:dist-p_j}
    p_{J}[v] = \frac{|\{i\in J: v\in X_i\}|}{|J| r} \qquad \forall v\in V
\end{equation}
(If all subsets in $\cX$ are vertex-disjoint, then $ p_{J}[v]$ is simply the uniform distribution over vertices in $\cup_{i\in I_t}X_i$ conditional on $I_t = J$; the distribution in \Cref{eq:dist-p_j} is simply adjusted for the overlap between the subsets in $\cX_{I_t}$. Intuitively, we select a uniform random set $X_i \in \cX_{I_t}$, and then select a uniform random vertex in $X_i$.) In the next step, we want the walk $W$ to get to a random vertex $v_{t+1}$ in a set in $\cX_{\Bar{I}_t}$, whose indices have not been selected; $v_{t+1}$ is  distributed according to:

\begin{equation}\label{eq:dist-p_j+1}
    q_{J}[v] = \frac{|\{i\in \Bar{J}: v\in X_i\}|}{|\Bar{J}| r} \qquad \forall v\in V
\end{equation}
which is also an adjusted uniform distribution to account for the overlap between the subsets in $\cX_{\Bar{I}_t}$. By moving from $v_t$ to $v_{t+1}$ along the shortest path $\pi_G(v_t, v_{t+1})$, $W$ hits a newly selected subset; see \Cref{fig:wasserstein-walk}.
If these distributional properties were all we needed,
then the algorithm and its analysis could be quite simple.
However, we simultaneously need a bound on the expected total length of the walk,
which we will accomplish by ensuring each component of the walk is short.
In particular, our goal here is to move from $v_t$ to $v_{t+1}$ (sampling from their respective distributions in a correlated manner), so that $\bbE[|\pi_G(v_t,v_{t+1})|]$ is small. The minimum possible expectation is  \emph{exactly} the Wasserstein distance, defined formally in \Cref{def:wasserstein}, between $p_{I_t}[v]$ and $q_{I_t}[v]$. (This problem has been studied extensively in optimal transport literature, but our construction does not require any deep understanding of this area.) As we will see, the solution to this problem gives us a distribution $y^*_t$ on $V$, from which we can sample  $v_{t+1}$. Finally, we walk from $v_t$ to $v_{t+1}$ and proceed to the next iteration. The pseudocode is given in \Cref{alg:WassertainWalk}.

\begin{definition}[Wasserstein distance]
 \label{def:wasserstein}
For a finite metric space $(U,d_U)$, the \defn{Wasserstein distance} between two probability distributions $p,q:U\to[0,1]$ (so $\sum_{u\in U}p[u]=\sum_{u\in U}q[u]=1$) with respect to $d_U$, denoted by \EMPH{$W_d(p,q)$}, is the solution of the following linear program.
\begin{equation}\label{eq:Wass-LP}
W_d(p,q) \coloneqq
\begin{array}[t]{@{}r@{\quad}r@{}l@{\quad}l@{}}
\min\limits_{x\in\real^{U\times U}\mathstrut} & \multicolumn{3}{@{}l@{}}{\sum_{u,v\in U}d_U(u,v)\,x[u,v]}\\
\textup{s.t.} & \sum_{v\in U}x[u,v] & {}=p[u] & \forall u\in U\\
& \sum_{u\in U}x[u,v] & {}=q[v] & \forall v\in U\\
& x & {}\geq0
\end{array}
\end{equation}
\end{definition}

A simple interpretation of the Wasserstein distance in \Cref{def:wasserstein} is this: fix a vertex $u\in U$; if we move $u$ to $v$ according to the distribution:
\begin{equation*}
    \frac{x[u,v]}{\sum_{v'\in U}x[u,v']} = \frac{x[u,v]}{p[u]}
\end{equation*}
then the expected cost of moving from $u$ to a random $v$ is $\sum_{v\in U}d_U(u,v)\frac{x[u,v]}{p[u]}$. If we move from a random $u$ sampled from $p$ to a random $v$ sampled from $q$, then the expected cost is:
\begin{equation}
    \sum_{u\in U} \left(\sum_{v\in U} d_U(u,v)\frac{x[u,v]}{p[u]}\right)p[u] = \sum_{u,v\in U}d_U(u,v)\,x[u,v] 
\end{equation}

That is, finding $W_d(p,q)$ is equivalent to minimizing the expected distance from a random $u$ to a random $v$, distributed according to $p$ and $q$, respectively.

\begin{algorithm}[H]
\caption{\textsc{WassersteinWalk}($G,\cX \coloneqq\{X_1,\ldots, X_k\}$)}
\label{alg:WassertainWalk}
\begin{algorithmic}[1]
\Require a connected graph $G=(V,E)$ and subsets $X_1,\dots,X_k\subset V$ of size $r = \ceil{n/k}$.
\Ensure a walk $W$ through $G$ that intersects every subset in $\cX$.
\State sample $i\in [k]$ uniformly, and sample $v_1\in X_i$ uniformly.
\State $I_1\gets\{i\}$ and $W\gets \{v_1\}$. 
\For{$t\leftarrow 1$ to $k-1$} \Comment{\lightgray{$v_t$ current endpoint of $W$ and $I_{t}$ indices visited up to step $t$}\label{line:forloop}}
    \State $p_{I_t}[v] \gets |\{i\in I_t:v\in X_i\}|/(|I_t|\cdot r)$  for all $v\in V$ \Comment{\lightgray{distribution in \Cref{eq:dist-p_j}}\label{line:pj}}
    \State  $q_{I_t}[v] \gets |\{i\in \Bar{I}_t:v\in X_i\}|/(|\Bar{I}_t|\cdot r)$ for all $v\in V$\Comment{\lightgray{distribution in \Cref{eq:dist-p_j+1}}}\label{line:pj+1}
    \State $x^*_t\leftarrow  $ optimal solution of \Cref{eq:Wass-LP} for computing Wasserstein distance $W_{d_G}(p_{I_t},q_{I_t})$
    \State $y_{t}^*[v] \gets
     \frac{x^*_t[v_{t},v]}{\sum_{v'\in V} x^*_t[v_{t},v']}
$ $\quad$ for all $v\in V$ \Comment{\lightgray{the distribution to sample $v_{t+1}$}}
    \State sample $v_{t+1}\in V$ according to $y^*_t$
    \State sample $i_t \in \{i\in\Bar{I}_t: v_{t+1} \in X_i\}$ uniformly \Comment{\lightgray{The walk $W$ will hit $X_{i_t}$}} 
    \State $I_{t+1}\leftarrow I_t \cup \{i_t\}$
    \State $W\leftarrow W\circ \pi_G(v_{t}, v_{t+1})$ 
\EndFor
\State \Return $W$
\end{algorithmic}
\end{algorithm}

In \Cref{subsec:Wass-vs-expansion}, we give a formal bound on the Wasserstein distance using the vertex expansion of graphs. Using this bound, we can bound the expected length of $\pi_G(v_t,v_{t+1})$ at the $t$-th iteration, as stated in the following lemma.

\begin{restatable}{lemma}{steplength}
\label{lemma:step-length}
At iteration $t\in[k-1]$ of \Cref{alg:WassertainWalk},
\begin{equation}
    \bbE
    \bigl[|\pi_G(v_t,v_{t+1})|\bigr]
    \leq \frac{4k}{\psi(G)\sqrt{t(k-t)}}.
\end{equation}
Here the expectation is over the random choices of the algorithm.
\end{restatable}

Since proving \Cref{lemma:step-length} requires understanding the relationship between Wasserstein distance and vertex expansion laid out in \Cref{subsec:Wass-vs-expansion}, we defer it to \Cref{subsec:expected-dist}. Here, we show that \Cref{lemma:step-length} implies \Cref{thm:hitting}.

\begin{proofof}{\Cref{thm:hitting}}
Assume without loss of generality (by discarding elements) that
each of the sets $X_i$ has exactly $\lceil n/k\rceil$ vertices. Let $W$ be the walk from running \Cref{alg:WassertainWalk} on $G$ and $\{X_1,\ldots, X_k\}$. By \Cref{lemma:step-length}:

\begin{equation}
\begin{split}
    \bbE[|W|] &\leq     \frac{4k}{\psi(G)}
\sum_{t=1}^{k-1}\frac1{\sqrt{t(k-t)}} \quad \text{(linearity of expectation)}\\
            &\leq \frac{4k}{\psi(G)}
\int_0^k\frac{dx}{\sqrt{x(k-x)}} \\
&= \frac{4\pi k}{\psi(G)} \qquad \text{(by the calculation in \Cref{eq:integral} below)}\\
&< \frac{13k}{\psi(G)}.
\end{split}
\end{equation}

The second inequality is because $f(x) \coloneqq 1/\sqrt{x(k-x)}$ decreases on $(0,k/2]$ and increases on $[k/2,k)$.
It remains to show that $\int_0^k\frac{dx}{\sqrt{x(k-x)}}=\pi$.
This follows by substituting $x:=k\sin^2\theta$ where $\theta$ ranges from $0$ to $\pi/2$. Since $dx=2k\cdot\sin\theta\cdot\cos\theta\cdot d\theta$ and $\sqrt{x(k-x)}=\sqrt{k\sin^2\theta\cdot k\cos^2\theta}=k\sin\theta\cdot\cos\theta$, we have
\begin{equation}\label{eq:integral}
    \int_0^k\frac{dx}{\sqrt{x(k-x)}}
=\int_0^{\pi/2}\frac{2k\sin\theta\cdot\cos\theta}{k\sin\theta\cdot\cos\theta}\,d\theta
=\int_0^{\pi/2}2\,d\theta
=\pi,
\end{equation}
as claimed. Finally, because $\E[|W|] < \frac{13k}{\psi(G)}$, there must exist some execution of \Cref{alg:WassertainWalk} that produces a walk with length less than $\frac{13k}{\psi(G)}$.
\end{proofof}

\subsection{Wasserstein distance and vertex expansion}\label{subsec:Wass-vs-expansion}

This section proves the following lemma, which bounds the Wasserstein distance between $p$ and $q$ in terms of $\|p-q\|_2$ and the graph expansion. This lemma is the only fact needed to prove \Cref{lemma:step-length}, so readers who are interested in the proof of \Cref{lemma:step-length} can skip this section.

\begin{lemma}
\label{lemma:wasserstein}
Let $G=(V,E)$ be a connected graph on $n$ vertices with shortest-path metric $d_G$.
If $p,q:V\to[0,1]$ are probability distributions,
then
\begin{equation}
  W_{d_G}(p,q)\leq\frac{4\sqrt n\cdot\|p-q\|_2}{\psi(G)}.   
\end{equation}
\end{lemma}

Although we are not aware of the exact form of \Cref{lemma:wasserstein} in the literature, it can be easily deduced from two known results.
The first is an immediate consequence of
the well-known Kantorovich–Rubinstein dual of Wasserstein distance.
The Kantorovich–Rubinstein dual is itself equivalent to the routine dual linear program for $W_{d_G}(p,q)$
in \Cref{eq:Wass-LP}
so long as the weights $d_U(u,v)$ satisfy the triangle inequality,
and it can be stated in the following form for finite metrics.

\begin{lemma}[{Kantorovich–Rubinstein dual, alternative form~\cite[Equation~(6.2)]{gabriel2019computational}}]\label{lm:Wass-dual}
   \begin{equation}
       W_d(p,q) = \begin{array}[t]{@{}r@{\quad}r@{}l@{\quad}l@{}}
\max\limits_{f\in\real^U\mathstrut} & \multicolumn{3}{@{}l@{}}{\sum_{u\in U}(p[u]-q[u])f_u}\\
\textup{s.t.} & |f_u-f_v| & {}\leq d_U(u,v) & \forall u,v\in U.
\end{array}
   \end{equation}
\end{lemma}

The second result is by  Bobkov, Houdr{\'{e}} and Tetali~\cite{BHT00},
relating
vertex expansion to feasible solutions of the dual program in \Cref{lm:Wass-dual} applied to shortest-path metrics.

\begin{theorem}[{\cite[Theorem~1]{BHT00}}]
\label{thm:poincare}
Let $G=(V,E)$ be a connected graph on $n$ vertices with shortest-path metric $d_G$.
If $f\in\real^V$ satisfies $|f_u-f_v|\leq d_G(u,v)$ for all $u,v\in V$, then
\[
\sum_{v\in V}(f_v-\overline f)^2\leq\frac{2n}{\left(\sqrt{1+\psi(G)}-1\right)^2},
\qquad\text{where }\overline f:=\frac1n\sum_{v\in V}f_v.
\]
\end{theorem}

\begin{proofof}{\Cref{lemma:wasserstein}}
Let $f$ be an optimal solution of the dual program in \Cref{lm:Wass-dual} where the metric is $(V,d_G)$. We have $W_{d_G}(p,q)=\sum_{v\in V}(p_v-q_v)f_v$ and $|f_u-f_v|\leq d_G(u,v)$ for every $u,v\in V$.  By \Cref{thm:poincare}:
\begin{equation} \label{eq:sum-square}
    \begin{split}
        \sum_{v\in V}(f_v-\overline f)^2
&\leq\frac{2n}{\left(\sqrt{1+\psi(G)}-1\right)^2}\\
&=\frac{2n\cdot\left(\sqrt{1+\psi(G)}+1\right)^2}{\left(\left(\sqrt{1+\psi(G)}-1\right)\left(\sqrt{1+\psi(G)}+1\right)\right)^2}\\
&=\frac{2n\cdot\left(\sqrt{1+\psi(G)}+1\right)^2}{\left(\left(1+\psi(G)\right)-1\right)^2}\\
&\leq\frac{2n\cdot\left(\sqrt3+1\right)^2}{\psi(G)^2} \qquad \text{(since $\psi(G)\leq 2$)}\\
&<\frac{16n}{\psi(G)^2}.
    \end{split}
\end{equation}

Since $\sum_{v\in V}p[v] =\sum_{v\in V}q[v]=1$,
we have $\sum_{v\in V}(p[v]-q[v])\cdot\overline f=0$.
Therefore,
\begin{alignat*}{2}
        W_{d_G}(p,q)
 &=\sum_{v\in V}(p[v]-q[v])f_v \\
 &=\sum_{v\in V}(p[v]-q[v])(f_v-\overline f) \\
&\leq\|p-q\|_2\cdot\left(\sum_{v\in V}(f_v-\overline f)^2\right)^{1/2} &\qquad& \text{(by Cauchy-Schwarz)}\\
&\leq\frac{4\sqrt n\cdot\|p-q\|_2}{\psi(G)} &\qquad& \text{(by \Cref{eq:sum-square})}
\end{alignat*}
as claimed.
\end{proofof}

\subsection{Bounding the expected distance}\label{subsec:expected-dist}

In this subsection, we prove \Cref{lemma:step-length}. First, we show that \(I_t\) is uniformly distributed over the size-\(t\) subsets of \([k]\), and that conditional on \(I_t = J\), the endpoint \(v_t\) has distribution \(p_J\). All the probabilities and expectations are over the random choices of the algorithm.

\begin{lemma}
\label{lemma:distribution} For every $t\in [k-1]$ and every size-$t$-subset $J\subset[k]$, and every vertex $a\in V$,
\[
\Pr[I_t=J]=\frac1{\binom k t}
\qquad\text{and}\qquad
\Pr[v_t=a\mid I_t=J]=\frac{|\{i\in J:a\in X_i\}|}{t\cdot r}.
\]
\end{lemma}
\begin{proof}

We proceed by induction on $t$.  The base case $t=1$ holds by lines 1 and 2 of \Cref{alg:WassertainWalk}.  Now assume that $t<k-1$ and that the claim holds for $t$.
By induction hypothesis, $v_t$ is distributed according to $p_J$. \Cref{alg:WassertainWalk} constructs a coupling between the distributions $p_J$ and $q_J$ and uses this to sample $v_{t+1}$, so $v_{t+1}$ has distribution $q_J$: formally,
\[\Pr[v_{t+1} = b | I_t = J] = \sum_{a \in V} x_t^*[a,b] = q_J[b] \qquad \text{for all $b \in V$}.\]
Now, notice that (conditional on $I_t = J$) the algorithm first samples $v_{t+1}$ according to $q_J$, and then samples $i_t$ uniformly at random among $\set{i \in [k] \setminus J : v_{t+1} \in X_i}$. This leads to a key observation: for any $i \in [k] \setminus J$ and $b \in X_i$,
\begin{equation}
\label{eq:joint-dist}
  \Pr[i_t = i \land v_{t+1} = b | I_t = J] = q_J[b] \cdot \frac{1}{|\set{i \in [k] \setminus J : b \in X_i}|} = \frac{1}{(k-t)r}.  
\end{equation}

\paragraph{Part 1 of the claim.} We prove that $I_{t+1}$ is uniformly distributed. First observe that (conditional on $I_t = J$), the next selected index $i_t$ is uniformly distributed among $[k] \setminus J$. Indeed, we can sum \eqref{eq:joint-dist} over the $r$ vertices of $X_i$:
\[\Pr[i_t = i | I_t = J] = \sum_{b \in X_i} \Pr[i_t = i \land v_{t+1} = b | I_t = J] = \frac{1}{k-t}.\]
The first part of the claim now follows: if we start with $I_t$ as a uniform random $t$-element subset of $[k]$ (which we do, by induction), and add a uniformly random element in $[k] \setminus I_t$ to form $I_{t+1}$, the set $I_{t+1}$ is a uniformly random $(t+1)$-element subset of $[k]$. That is, $\Pr[I_{t+1} = J'] = \frac{1}{\binom{k}{t+1}}$ for any size-$(t+1)$ subset $J' \subseteq [k]$.\footnote{In more detail: the event $I_{t+1} = J'$ occurs precisely when $I_t = J' \setminus \set{i}$ and $i_t = i$ for some $i \in J'$. For each $i \in J'$, this event occurs with probability $\frac{1}{(k-t)\binom{k}{t}}$. Summing over the $t+1$ choices of $i$ yields $\Pr[I_{t+1} = J'] = \frac{1}{\binom{k}{t+1}}$}

\paragraph{Part 2 of the claim.} We prove that the distribution of $v_{t+1}$ has the desired form.  
Let $J'\subset[k]$ be a size-$(t+1)$ subset of $[k]$ and let $a\in V$.
Observe that
\begin{align*}
\Pr[I_{t+1}=J'\land v_{t+1}=a]
&=\sum_{i\in J':a\in X_i}\Pr[I_t=J'\setminus\{i\}\land v_{t+1}=a\land i_t=i]\\
&=\sum_{i\in J':a\in X_i}\frac1{\binom kt}\cdot\frac1{(k-t)\cdot r}
=\frac{|\{i\in J':a\in X_i\}|}{\binom kt\cdot(k-t)\cdot r}
\end{align*}
where the second-to-last equality follows from \eqref{eq:joint-dist}. Hence
\begin{align*}
\Pr[v_{t+1}=a\mid I_{t+1}=J']
&=\frac{\Pr[I_{t+1}=J'\land v_{t+1}=a]}{\Pr[I_{t+1}=J']}\\
&=\frac{|\{i\in J':a\in X_i\}|}{\binom kt\cdot(k-t)\cdot r}\cdot {k \choose t+1}\\
&=\frac{|\{i\in J':a\in X_i\}|}{(t+1)\cdot r}~
\end{align*}
as claimed.
\end{proof}

Now we prove \Cref{lemma:step-length} restated below.

\steplength*
\begin{proof} By \Cref{lemma:distribution}, given $I_t = J$, we have $v_{t}=a$ with probability $p_J[a]$. Then in line 8, we sample $v_{t+1}$ and:
\begin{equation}
    \pr[v_{t+1} = b\mid  I_t = J \wedge v_t = a] = \frac{x^*_t[a,b]}{\sum_v x^*_t[a,v]} = \frac{x^*_t[a,b]}{p_J[a]}.
\end{equation}
Then, we have:
\begin{equation}
\label{eq:conditional}
\begin{split}
    \bbE[d_G(v_t,v_{t+1})\mid I_t = J]
&=\sum_{a,b\in V}d_G(a,b)\cdot p_J[a]\cdot\frac{x^*_t[a,b]}{p_J[a]}\\
& =W_{d_G}(p_J,q_{J})\\
&\leq\frac{4\sqrt n\cdot\|p_J -q_{J}\|_2}{\psi(G)} \qquad \text{(by \Cref{lemma:wasserstein})}
\end{split}
\end{equation}
By averaging over $I_t$ (rather than conditioning on $I_t = J$), we get:
\begin{equation}\label{eq:expected-dist-1}
\begin{split}
    \bbE\left[d_G(v_t,v_{t+1})\right]
&\leq\frac{4\sqrt n}{\psi(G)}\cdot\bbE\left[\|p_{I_t}-q_{I_t}\|_2\right]\\
&\leq\frac{4\sqrt n}{\psi(G)}\cdot\sqrt{\bbE\left[\|p_{I_t}-q_{I_t}\|_2^2\right]}
\qquad \text{(by Jensen's inequality)}    
\end{split}
\end{equation}
Now we focus on bounding $\bbE\left[\|p_{I_t}-q_{I_t}\|_2^2\right]$. To this end, for each $i\in[k]$, let $\sigma_i:=1/t$ if $i\in I_t$ and $\sigma_i:=-1/(k-t)$ otherwise. Then by the definition of $p_{I_t}$ and $q_{I_t}$ in \Cref{eq:dist-p_j,eq:dist-p_j+1}, for every  $v\in V$,
\begin{equation}
\label{eq:difference}
p_{I_t}[v]-q_{I_t}[v] =\frac1r\sum_{i:v\in X_i}\sigma_i.
\end{equation}
By \Cref{lemma:distribution}, $I_t$ is a uniformly random size-$t$ subset of $[k]$,
so each $i\in[k]$ lies in $I_t$ with probability $t/k$, and
\begin{equation}\label{eq:sigma-sq}
    \bbE[\sigma_i^2]
=\frac tk\cdot\frac1{t^2}+\frac{k-t}k\cdot\frac1{(k-t)^2}
=\frac 1{tk}+\frac{1}{k(k-t)}
=\frac1{t(k-t)}.
\end{equation}
For distinct $i,i'\in[k]$, considering whether both, neither or exactly one of $i$ and $i'$ lie in $I_t$ gives
\begin{equation} \label{eq:negative-corr}
    \begin{split}
        \bbE[\sigma_i\sigma_{i'}]
&=\frac{t(t-1)}{k(k-1)}\cdot\frac1{t^2}
+\frac{(k-t)(k-t-1)}{k(k-1)}\cdot\frac1{(k-t)^2}
-\frac{2t(k-t)}{k(k-1)}\cdot\frac1{t(k-t)}\\
&=\frac1{k(k-1)}\left(\frac{t-1}t+\frac{k-t-1}{k-t}-2\right)
=-\frac1{k(k-1)}\left(\frac1t+\frac1{k-t}\right)
<0.
    \end{split}
\end{equation}

By \Cref{eq:difference}, we get
\begin{equation}\label{eq:expected-norm-sq}
    \begin{split}
        \bbE\left[\|p_{I_t}-q_{I_t}\|_2^2\right]
&=\frac1{r^2}\sum_{v\in V}\bbE\left[\left(\sum_{i:v\in X_i}\sigma_i\right)^2\right]
=\frac1{r^2}\sum_{v\in V}\sum_{i,i':v\in X_i\cap X_{i'}}\bbE[\sigma_i\sigma_{i'}]
\\
&\leq\frac1{r^2}\sum_{v\in V}\sum_{i:v\in X_i}\bbE[\sigma_i^2] \qquad \text{(by \Cref{eq:negative-corr})}
\\
&=\frac1{r^2}\sum_{i\in[k]}\bbE[\sigma_i^2]\cdot|X_i|
\\
&=\frac1{r^2}\cdot k\cdot\frac1{t(k-t)}\cdot r \qquad \text{(by \Cref{eq:sigma-sq})}\\
&=\frac k{t(k-t)\cdot r}.
\end{split}
\end{equation}

Plugging \Cref{eq:expected-norm-sq} into \Cref{eq:expected-dist-1}, we get:
\begin{align*}
\bbE\left[d_G(v_t,v_{t+1})\right]
&\leq\frac{4\sqrt n}{\psi(G)}\cdot\sqrt{\frac k{t(k-t)\cdot r}}
&&\\
&\leq\frac{4\sqrt{k\cdot r}}{\psi(G)}\cdot\sqrt{\frac k{t(k-t)\cdot r}}
&&\text{since $r\geq n/k$}\\
&=\frac{4k}{\psi(G)\sqrt{t(k-t)}}.
\end{align*}
as claimed.
\end{proof}

\section{The Algorithm}

For a subgraph $H$ of $G$ and a set of vertices $A$ of $G$, let \EMPH{$\bdry_H A$} denote the set of vertices of $H$ that are not in $A$ but have at least one neighbor in $A$. We note that $A$ might not belong to $V(H)$; for example, when $H$ is a subgraph of a larger graph $G$ and $A\subseteq V(G)$ that is not in $H$. When $A \subseteq V(H)$, we use the notation \EMPH{$B_H(A,r)$} to denote the ball of radius $r$ around $A$ in $H$, that is, $\set{v \in V(H) : \dist_H(A, v) \le r}$. We will use the following corollary of \Cref{thm:hitting}.

\begin{corollary}\label{cor:connector} Suppose that $H$ is a $\psi$-vertex expander.  Let $\cX = \{X_1,X_2,\ldots, X_k\}$  be a collection of $k$ vertex subsets of $V(H)$ such that $|X_i|\geq |V(H)|/h$ and $k\leq h$. Then there is a procedure \EMPH{$\textsc{ShortTree}(H,\cX)$} that produces a tree  $T\subseteq H$ such that:
\begin{enumerate}
    \item $T\cap X_i \not=\varnothing$ for every $i\in [k]$.
    \item  $|V(T)| \leq \Ceil{13h/\psi}$.
\end{enumerate}
\end{corollary}
\begin{proof}
    If $k = h$, then  we can apply \Cref{thm:hitting} to produce a walk $W$, and we can then simply take $T$ to be a spanning tree of $H[V(W)]$. As $W$ has length strictly less than $13 h /\psi$, the tree $T$ has at most $\Ceil{13h/\psi}$ vertices. If $k < h$, then we duplicate $X_k$ to get $\cX^+ = \{X_1,\ldots, X_{k}, X_{k+1} = X_k, 
    \ldots, X_h = X_k\}$ and apply \Cref{thm:hitting} to $H$ and $\cX^+$.
\end{proof}

\paragraph{The algorithm.~} We begin with an informal description of the algorithm. The formal pseudocode is given as \Cref{alg:separator}. At every step, the algorithm maintains three important objects: $H$, $X$ and $\cK$.

\begin{enumerate}
    \item \EMPH{$H$} is an induced subgraph of $G$. Except when the algorithm terminates, $H$ contains at least $3n/4$ vertices.  One should think of the ``current separator'' as \EMPH{$\bdry_G H$}, and $H$ as the largest connected component\footnote{For easier bookkeeping, $H$ may briefly become disconnected at the end of an iteration; it is reassigned to the largest connected component at the start of the next iteration.} of $G$ if we remove $\bdry_G H$ from $G$.
\end{enumerate}
The  ``current separator'' $\bdry_G H$ will lie in the union of two objects:
\begin{enumerate}
    \item[2.]  \EMPH{$X$} is a subset of vertices of $G$ which the algorithm will keep adding to. We guarantee that $|X| \le O(h \sqrt{n})$.
    \item[3.]  \EMPH{$\cK$} is a clique minor model found so far, whose branch sets of $\cK$ are disjoint from $H$. We guarantee that, for every $C \in \cK$, we have $|\bdry_G H \cap C| \le O(\sqrt{n})$.
\end{enumerate}

\begin{wrapfigure}{r}{0.35\textwidth}
    \centering
    \vspace{-10pt}
    \includegraphics[width=0.8\linewidth]{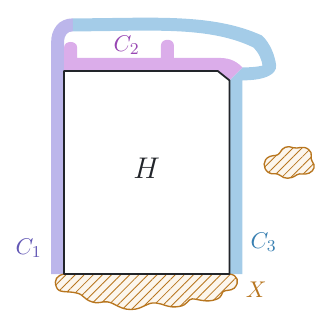}
    \caption{The subgraph $H$. The boundary of $H$ is contained in the union of $X$ (drawn in orange) and the vertices of the $K_3$ minor model $\set{C_1, C_2, C_3} = \cK$.}
    \label{fig:sep-1}
    \vspace{-10pt}
\end{wrapfigure}

It will hold that $X$ and $V(\cK)$ are disjoint from $V(H)$. See \Cref{fig:sep-1}. In every step, the algorithm removes at least one vertex from $H$ and at most $n/2$ vertices from $H$, adding them to $X$ or $V(\cK)$. The algorithm terminates when the size of $H$ drops below $3n/4$, at which point we return $\bdry_G H$ as the separator.\footnote{This description is slightly inaccurate, as $H$ may be disconnected at the start of an iteration. The first step of each iteration reassigns $H$ to be the largest connected component, and the termination condition is that every connected component of $H$ has $< 3n/4$ vertices. } Because each step removes at most $n/2$ vertices, we have $n/4 \le |H| \le 3/4 n$ and thus $\bdry_G H$ is a balanced separator. To guarantee that $\bdry_G H$ has size $O(h \sqrt{n})$, we maintain the invariant that $\bdry_G H \subseteq X \cup V(\cK)$, that $|X| \le O(h \sqrt n)$, and that for each branch set $C \in \cK$ we have $|C \cap \bdry_G H| \le O(\sqrt n)$ (recall that $K_h$-minor-freeness implies that $\cK$ contains at most $h-1$ branch sets).

\begin{figure}[htbp]
  \centering
  \begin{subfigure}[t]{0.3\textwidth}
    \centering
    \includegraphics[width=\linewidth]{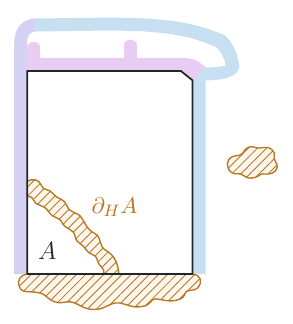}
    \caption{Case 1: Add $\partial_H A$ to $X$.}
    \label{fig:case-1}
  \end{subfigure}
  \hfill
    \textcolor{black!25}{\rule[1cm]{0.4pt}{3cm}}
    \hfill
  \begin{subfigure}[t]{0.3\textwidth}
    \centering
    \includegraphics[width=\linewidth]{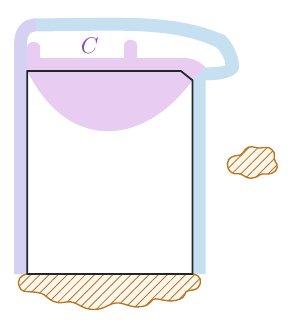}
    \caption{Case 2: Grow $C$ by
      $C\leftarrow C\cup B_H(\partial_H C,r)$.}
    \label{fig:case-2}
  \end{subfigure}
  \hfill
  \textcolor{black!25}{\rule[1cm]{0.4pt}{3cm}}
  \hfill
  \begin{subfigure}[t]{0.3\textwidth}
    \centering
    \includegraphics[width=\linewidth]{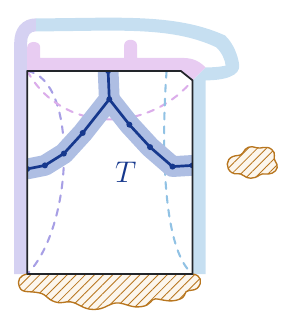}
    \caption{Case 3: Add a new branch set $T$ to $\mathcal K$ using the expander walk theorem.}
    \label{fig:case-3}
  \end{subfigure}
  \caption{The three cases possible in each step of the algorithm}
  \label{fig:sep-2}
\end{figure}

Each iteration of the algorithm proceeds as follows. At the start, the algorithm checks whether any connected component of $H$ has $\ge 3n/4$ vertices. If not, we terminate and return the current separator $\bdry_G H$; otherwise we update $H$ to be the largest connected component, update $\cK$ by discarding all branch sets that aren't adjacent to $H$, and proceed. Within each iteration, there are three cases; see \Cref{fig:sep-2}.
\begin{itemize}
    \item \textbf{Case 1: $H$ is not a $h/\sqrt{n}$-vertex expander.}
    
    By the definition of a vertex expander, there must be a subset of vertices $A$ such that $|\bdry_H A|\leq h/\sqrt{n} |A|$ and $|A| \leq |V(H)|/2$ (the smaller side of the sparse cut). We will update $X$ by adding $\bdry_H A$ to $X$, and removing $\bdry_H A$ from $H$. We can afford adding   $\bdry_H A$ to $X$ since we can charge the size of $\bdry_H A$ to $A$. Each vertex gets $h/\sqrt{n}$ charge over the entire process, since $A$ will be removed from $H$ (once $H$ is reassigned to the largest connected component) and so over the entire course of the algorithm, the size of $X$ is bounded by $(h/\sqrt{n})n = h\sqrt{n}$.
    \item  \textbf{Case 2:  There exists branch set $C \in \cK$ and $r \geq 1$ with $|B_H(\bdry_H C, r)|\leq n/2$ and $|\bdry_H B_H(\bdry_H C, r-1)| \leq \sqrt{n}$.}
    
    Then we simply grow $C$ by adding the entire ball $B_H(\bdry_H C, r)$ to $C$, and (correspondingly) deleting $B_H(\bdry_H C, r)$ from $H$. This step could blow up the size of $C$ significantly (up to $\Omega(n)$). However, since $|\bdry_H B_H(\bdry_H C, r-1)| \leq \sqrt{n}$, after the update we still have that $|C \cap \bdry_G H |\leq \sqrt{n}$. The condition $|B_H(\bdry_H C, r)|\leq n/2$
    just guarantees that we remove at most $n/2$ vertices from $H$ in this step.
    
    \item  \textbf{Case  3: Otherwise}. 
    
    In this case, (i) $H$ is an $h/\sqrt n$-vertex expander, and (ii) for every branch set $C \in \cK$ and radius $r\geq 1$, either $|B_H(\bdry_H C, r)|> n/2$ or $|\bdry_H B_H(\bdry_H C, r-1)| > \sqrt{n}$.  If we pick $r^* = \ceil{\sqrt{n}/h}$, we can show (\Cref{clm:n/h-size}) that $|B_H(\bdry_H C, r^*)| > n/h$. We can now form $\cX \leftarrow \{B_H(\bdry_H C, r^*)\}_{C\in \cK}$ and apply \Cref{cor:connector} to get a tree \EMPH{$T$} of $H$ of size $O(\sqrt{n})$ that intersects every set in $\cX$. We then \emph{extend} $T$ so that it is adjacent to every branch set in $\cK$: for every branch set $C\in\cK$,  we pick an arbitrary vertex $x_C\in T\cap B_H(\bdry_H C, r^*)$, and add to $T$ the shortest path in $H$ from $x_C$ to a vertex in $\bdry_H C$.  Note that since $r^* = \ceil{\sqrt{n}/h}$ and $|\cK|\leq h$, the total size of all the paths added to $T$ is at most $hr^* = O(\sqrt{n})$, so the total size of (a new branch set) $T$ formed in this step is still $O(\sqrt{n})$. We add $T$ to $\cK$, and remove $T$ from $H$.
\end{itemize}

In summary, the algorithm guarantees that $|X| = O(h\sqrt{n})$ (Case 1) and for every branch set $C\in \cK$, either $|\bdry_G H \cap C|\leq \sqrt{n}$ (Case  2) or directly that  $|C| = O(\sqrt{n})$  (Case 3).  Moreover, the current separator $\bdry_G H$ is in $X \cup V(\cK)$. Since $|\cK| < h$, the size of the current separator is always $O(h\sqrt{n})$. Once $|V(H)|\leq 3n/4$, the current separator is balanced and we simply return it.

\begin{algorithm}[H]
\caption{$\textsc{BalancedSeparator}(G)$\Comment{\lightgray{Find a balanced separator of size $O(h\sqrt{n})$}}}\label{alg:separator}
\begin{algorithmic}[1]
    \State  $\mathcal{K}\gets \varnothing, \,\, X\leftarrow \varnothing, \,\, H \leftarrow G$ 
    \While{there is some connected component of $H$ with $\ge 3n/4$ vertices}
        \State  $H\leftarrow$ largest connected component of $H$\label{line:connected}
        \State $\mathcal{K}\leftarrow \set{C \in \cK: \text{$C$ has a neighbor in $H$}}$\label{line:trim}

        \If{$H$ \underline{is not} a $\frac{h}{\sqrt{n}}$-vertex expander}
        \Comment{\blue{Case 1}}
            \State $A\leftarrow$ nonempty subset of $V(H)$ such that $|A|\leq n/2$ and $|\bdry_H A|\leq \frac{h}{\sqrt{n}}|A|$
            \State $X\leftarrow X\cup \bdry_H A$. \label{line:add-bdry}
            \State $H\leftarrow H \setminus \bdry_H A$
        \ElsIf{there exists $C \in \cK$ and $r \geq 1$ with $|B_H(\bdry_H C, r)|\leq n/2$ and $|\bdry_H B_H(\bdry_H C, r - 1)| \leq \sqrt{n}$ } 
    \Comment{\blue{Case 2}}
                \State $B \gets B_H(\bdry_H C, r)$
                \State $C\leftarrow C\cup  B$. \label{line:c-expand}
                \State $H\leftarrow H \setminus B$
        \Else \Comment{\lightgray{$\forall C \in \cK, \forall r\geq 1$: $|B_H(\bdry_H C, r)|> n/2$ or $|\bdry_H B_H(\bdry_H C, r-1)| > \sqrt{n}$ }}\Comment{\blue{Case 3}}
        \If{$\cK = \varnothing$}
            \State $T \gets $ an arbitrary vertex in $V(H)$\label{line:singleton}
        \Else
            \State  $r^*\leftarrow \ceil{\sqrt{n}/h}$
            \State $\cX \leftarrow \{B_H(\bdry_H C, r^*)\}_{C\in \cK}$ \label{line:CX-size} \Comment{\lightgray{\Cref{clm:n/h-size} shows that $|B_H(\bdry_H C, r^*)|\geq n/h$}.}
            \State $T\leftarrow \textsc{ShortTree}(H,\cX)$ \label{line:new-branchset}
            \For{every $C\in \cK$}
                \State $x_C\leftarrow $ a vertex in $T\cap B_H(\bdry_H C, r^*)$ \State $\pi_{x_C}\leftarrow $ a shortest path in $H$ from $x_C$ to a vertex in $\bdry_H C$ \Comment{\lightgray{$\pi_{x_C}$ exists because $H$ is connected.}}
                \State $T\leftarrow T\cup \pi_{x_C}$ \label{line:extend-T}  \Comment{\lightgray{this makes $T$ adjacent to $C\in \cK$.}}
            \EndFor
        \EndIf
            \State $\cK \gets \cK \cup \set{T}$
            \State $H \gets H \setminus V(T)$
        \EndIf
    \EndWhile
    \State \Return 
    $\bdry_G H$
\end{algorithmic}
\end{algorithm}

We now turn to the analysis. First, we show that every set $B_H(\bdry_H C, r^*)\in \cX$ in line~\ref{line:CX-size} has size at least $n/h$. This justifies the call to $\textsc{ShortTree}(H,\cX)$ in line~\ref{line:new-branchset}.

\begin{claim}\label{clm:n/h-size} Let $C\in \cK$ be a branch set such that for every $r\geq 1$, either $|B_H(\bdry_H C, r)| > n/2$ or $|\bdry_H B_H(\bdry_H C, r-1)| > \sqrt{n}$. Then
\begin{equation}
    |B_H(\bdry_H C, r^*)|\geq n/h \ge |V(H)|/h
\end{equation}
where $r^* = \ceil{\sqrt{n}/h}$. Consequently, $(H, \cX)$ is a valid input to the \textsc{ShortTree} procedure on line~\ref{line:new-branchset}.
\end{claim}
\begin{proof}
    If there exists $r\in [0,r^*]$ such that $|B_H(\bdry_H C, r)|\geq n/h$, then the claim holds, as $B_H(\bdry_H C, r) \subseteq B_H(\bdry_H C, r^*)$. Otherwise, by assumption, for $r \in [0,r^*-1]$ we have that $|\bdry_H B_H(\bdry_H C, r)| > \sqrt{n}$.
    Observing that $B_H(\bdry_H C, r^*)$ contains $\bdry_H B_H(\bdry_H C, r)$ for every $r \in [0,r^*-1]$, we conclude
    \begin{equation*}
        |B_H(\bdry_H C, r^*) | \geq \sum_{r=0}^{r^*-1} |\bdry_H B_H(\bdry_H C, r)| > r^*\sqrt{n} \ge n/h
    \end{equation*}
    as claimed. Thus, when we call $\textsc{ShortTree}(H, \cX)$ in line~\ref{line:new-branchset}, it holds that every $X \in \cX$ has $|X| \ge |V(H)|/h$ and that $H$ is an $h/\sqrt n$-vertex expander, as required by the definition of $\textsc{ShortTree}$.
\end{proof}

This shows that the algorithm is well-defined. We now observe that it terminates.
\begin{lemma}
\label{lem:terminate}
The algorithm $\textsc{BalancedSeparator}(G)$ terminates.
\end{lemma}
\begin{proof}
    It suffices to show that every iteration of the while loop removes at least one vertex from $H$. Indeed, in Case 1, we remove some nonempty subset $\bdry_H A$ from $H$. In Case 2, we remove some ball $B_H(\bdry_H C, r)$ from $H$; this ball is nonempty because it contains $\bdry_H C$, which is nonempty because (due to the trimming in line 4) $C$ has some neighbor in $H$. In Case 3, we remove the nonempty set $T$ from $H$.
\end{proof}

We record two simple invariants. 
\begin{lemma}\label{lem:invariants} At the start of every iteration, the following properties hold:
\begin{enumerate}
    \item $\cK$ is a clique minor model in $G$ (of size $|\cK|$).  Thus $|\cK| < h$.
    \item  $V(H)$ and $V(\cK)$ are disjoint.
\end{enumerate}
    
\end{lemma}
\begin{proof}
In the beginning of the algorithm, we have $\cK = \varnothing$, $X = \varnothing$, and $H = G$, and so the invariants are satisfied trivially.

\medskip \noindent \textbf{1)} By induction, assume that $\cK$ is a model of a clique minor in $G$ at the start of the current iteration. Discarding some branch sets in line~\ref{line:trim} does not affect this property. In the rest of the iteration (lines 5--24), we may modify $\cK$ during one of the three cases of the algorithm. In Case 1, $\cK$ is unchanged. In Case 2, some $C \in \cK$ is expanded (line~\ref{line:c-expand}) to include additional vertices in $V(H)$. But $C$ remains disjoint from the other branch sets after the expansion (because, by induction, $H$ was disjoint from the other branch sets), $C$ remains adjacent to all other branch sets, and $C$ remains connected; thus $\cK$ remains a clique minor model. Finally, in Case 3, a new branch set $T$ is added. If $T$ is a singleton added by line~\ref{line:singleton} then the claim holds. In the nontrivial case (lines 17--23), the branch set $T$ is disjoint from the other $C \in \cK$ because $V(T) \subseteq V(H)$ and $V(H)$ is disjoint from the other branch sets; moreover, $T$ is adjacent to every other branch set $C \in \cK$ because we added the path $\pi_{x_C}$ to $T$ in line~\ref{line:extend-T}.

\medskip \noindent \textbf{2)} By induction, assume $V(\cK)$ and $V(H)$ are disjoint at the start of the current iteration. In this iteration, the set $V(H)$ shrinks, and some vertices may be added to $V(\cK)$. It suffices to show that every vertex added to $V(\cK)$ in the iteration is removed from $V(H)$. Indeed, in Case 1, no new vertices are added to $V(\cK)$; in Case 2, the ball $B = B_H(\bdry_H C, r)$ is added to $V(\cK)$ and removed from $V(H)$; and in Case 3, the vertices of the new branch set $T$ are added to $V(\cK)$ and removed from $V(H)$.
\end{proof}

We can now prove that the algorithm produces a separator of size $O(h \sqrt n)$ (\Cref{lem:size}), and that this separator is balanced (\Cref{lem:balance}). Our main theorem follows.

\begin{lemma}[Size bound]
\label{lem:size}
    At the start and end of every iteration of the while loop, we have $|\bdry_G H| \le 16 h \sqrt n$. This is a consequence of the following invariants:
    \begin{enumerate}
        \item $\bdry_G H \subseteq X \cup V(\cK)$
        \item $|X| \le h \sqrt n$
        \item For every $C \in \cK$, we have $|\bdry_G H \cap C| \le 15 \sqrt n$
    \end{enumerate}
\end{lemma}
\begin{proof}
Given the invariants, we can conclude that $|\bdry_G H| \le |X| + |\bdry_G H \cap V(\cK)| \le h \sqrt n + (h-1) \cdot 15 \sqrt n \le 16 h \sqrt n$, where the second-to-last inequality follows from the fact that $|\cK| \le h-1$ (by \Cref{lem:invariants}(1)). It remains to prove the invariants.

\medskip \noindent \textbf{1)} We proceed by induction. Assume that $\bdry_G H \subseteq X \cup V(\cK)$ holds at the start of the while loop (line 3). Observe that for any connected component $H'$ of $H$, we have $\bdry_G H' \subseteq \bdry_G H$; thus the claim holds after line~\ref{line:connected} (when the largest connected component of $H$ is renamed as $H$), and it also holds trivially after line~\ref{line:trim} (when all branch sets not adjacent to $H$ are discarded). In lines 5--25, in each of the three cases $H$ is modified in the following way: some set of vertices is added to $X \cup V(\cK)$, and these same vertices are removed from $H$. When removing a vertex $v$ from $H$, the only possible new vertex that could be added to $\bdry_G H$ is $v$ itself. Thus, after each iteration it remains the case that $\bdry_G H \subseteq X \cup V(\cK)$.

\medskip \noindent \textbf{2)} We prove $|X| \le h \sqrt n$ by a charging argument. Note that we only add to $X$ in line~\ref{line:add-bdry}: whenever we add $\bdry_H A$ to $X$, we charge the size of $\bdry_H A$ to $|A|$. Each vertex in $A$ gets at most $h/\sqrt{n}$ charge. Now we argue that vertices of $A \subset V(H)$ will not be charged again in future iterations. Indeed, after we delete $\bdry_H A$ from $H$ in this iteration, any connected component that intersects $A$ (and is therefore fully contained in $A$) contains at most $n/2$ vertices. But if the algorithm continues, it reassigns $H$ (in line~\ref{line:connected}) to be some connected component with at least $3n/4$ vertices. Thus vertices in $A$ are never charged again, and the total size of $X$ over all iterations is at most $(h/\sqrt{n})n = h\sqrt{n}$.

\medskip \noindent \textbf{3)} At the start of the algorithm, we have $\cK = \varnothing$ and the claim holds trivially.
There are three things that can occur during the course of the algorithm that could affect whether the claim holds: (1) the set $H$ can change, (2) some set $C \in \cK$ may be modified in Case 2, or (3) a new set $C \in \cK$ may be created in Case 3. We claim that after any of these changes, the claim still holds.

For the first, observe that $H$ only shrinks over the course of the algorithm (i.e., $H$ is changed to be some subgraph of itself). As every $C \in \cK$ is disjoint from $H$ (by \Cref{lem:invariants}(2)), shrinking $H$ can only cause the set $\bdry_G H \cap C$ to shrink.

 For the second, suppose that Case 2 selects some ball $B_r = B_H(\bdry_H C, r)$ and updates $C \gets C \cup B_r$ and $H \gets H \setminus B_r$. For clarity, let $H$ and $C$ denote the sets before the update, and let $H_{\rm new}$ and $C_{\rm new}$ denote the sets after the update. We want to show that $|\bdry_G H_{\rm new} \cap C_{\rm new}| \le \sqrt n$. 
 To this end, we first define $B_{r-1} = B_H(\bdry_H C, r-1)$ and observe that $B_{r} \setminus B_{r-1} = \bdry_H B_H(\bdry_H C, r-1)$.  We claim that
 \begin{equation}
 \label{eq:ring}
\bdry_G H_{\rm new} \cap C_{\rm new} \subseteq B_r \setminus B_{r-1}.     
 \end{equation}
The desired size bound on $\bdry_G H_{\rm new} \cap C_{\rm new}$ follows from \eqref{eq:ring}, as by assumption of Case 2 we have $|B_r \setminus B_{r-1}| \le \sqrt{n}$. To prove \eqref{eq:ring}, observe that $C_{\rm new} = C \cup B_r$, so it suffices to show that vertices in $C \cup B_{r-1}$ are not in $\bdry_G H_{\rm new}$. Indeed, every vertex in $C$ has no neighbors in $H_{\rm new}$, as every neighbor lying in $H$ is in $\bdry_H C \subseteq B_r$; similarly, every vertex in $B_{r-1}$ has no neighbors in $H_{\rm new}$, as its neighbors in $H$ are in $B_r$. So, the only vertices in $C_{\rm new}$ that could be adjacent to $H_{\rm new}$ (and thus in $\bdry_G H_{\rm new}$) are those in $B_r \setminus B_{r-1}$, as desired.
 
 Finally, for the third case, suppose that Case 3 adds some new branch set $T$ to $\cK$. When $T$ is formed in line~\ref{line:new-branchset}, by \Cref{cor:connector}
\begin{equation*}
    |V(T)| \leq \frac{13h}{h/\sqrt{n}} + 1 \le 13 \sqrt n + 1.
\end{equation*}    
Then $T$ is extended further in line \ref{line:extend-T} by adding at most $r^* = \ceil{\sqrt{n}/h}$ vertices per branch set, totaling at most $(h-1)\ceil{\sqrt{n}/h} \le \sqrt{n} - \sqrt{n}/h + h - 1 \leq 2 \sqrt{n} -1$. (If $h \ge \sqrt n$, our main theorem holds trivially.) Thus, $|V(T)|\leq 15\sqrt{n}$.
\end{proof}

\begin{lemma}[Balance]
\label{lem:balance}
    At every point during the algorithm, we have $|H| \ge n/4$. Moreover, when the algorithm terminates, every connected component of $H$ has at most $3n/4$ vertices.
\end{lemma}
\begin{proof}
    The second statement follows immediately from the while loop termination condition. We now show that $|H| \ge n/4$ always.
    Indeed, initially we have $|H| = n$. At the beginning of each iteration, after line~\ref{line:connected}, we have $|H| \ge 3n/4$ by assumption. During the iteration, in steps 4--25 of the algorithm, at most $n/2$ vertices are removed from $H$: in case 1 we remove $|\bdry_H A|$ vertices, where $|\bdry_H A| \le h/\sqrt n |A| \le |A| \le n/2$ (note we assume $h \le \sqrt n$ as otherwise our main theorem holds trivially); in Case 2 we remove $|B_H(\bdry_H C, r)| \le n/2$ vertices; and in Case 3 we remove $|V(T)|$ vertices, where the proof of \Cref{lem:size}(3) shows that $|V(T)| \le 15\sqrt{n}$, which is $ \le n/2$ if $n$ is sufficiently large\footnote{In the extreme case when $15 \sqrt n \ge n/2$, our main theorem holds trivially: deleting any set of $14 \sqrt n$ vertices is a balanced separator.}. Thus the size of $H$ is always at least $3n/4 - n/2 = n/4$.
\end{proof}

The main theorem follows immediately from \Cref{lem:size,lem:balance}.

\SepMain*
\begin{proof}
The algorithm returns $\bdry_G H$, which comprises at most $16 h \sqrt n$ vertices by \Cref{lem:size}. Moreover, every connected component of $G \setminus \bdry_G H$ is either a connected component of $H$ or is contained in $G \setminus H$; as \Cref{lem:balance} implies that connected components of $H$ have size $\le 3n/4$ and that $|V(G\setminus H)| \le 3n/4$, we conclude that $\bdry_G H$ is a balanced separator of $G$.
\end{proof}

\bibliographystyle{alphaurl}
\bibliography{ref}

\appendix
\end{document}